# EMPIRICAL RISK MINIMIZATION IN INVERSE PROBLEMS[1]

By Jussi Klemelä and Enno Mammen

*University of Oulu and University of Mannheim*

We study estimation of a multivariate function $f: \mathbf{R}^d \to \mathbf{R}$ when the observations are available from the function $Af$, where $A$ is a known linear operator. Both the Gaussian white noise model and density estimation are studied. We define an $L_2$-empirical risk functional which is used to define a $\delta$-net minimizer and a dense empirical risk minimizer. Upper bounds for the mean integrated squared error of the estimators are given. The upper bounds show how the difficulty of the estimation depends on the operator through the norm of the adjoint of the inverse of the operator and on the underlying function class through the entropy of the class. Corresponding lower bounds are also derived. As examples, we consider convolution operators and the Radon transform. In these examples, the estimators achieve the optimal rates of convergence. Furthermore, a new type of oracle inequality is given for inverse problems in additive models.

**1. Introduction.** We consider estimation of a function $f: \mathbf{R}^d \to \mathbf{R}$ when a linear transform $Af$ of the function is observed under stochastic noise. We consider both the Gaussian white noise model and density estimation with i.i.d. observations. We study two estimators: a $\delta$-net estimator which minimizes the $L_2$-empirical risk over a minimal $\delta$-net of a function class and a dense empirical risk minimizer which minimizes the empirical risk over the whole function class without restricting the minimization over a $\delta$-net. We call this estimator a "dense minimizer" because it is defined as a minimizer over a possibly uncountable function class. The $\delta$-net estimator is more universal: it may also be applied for nonsmooth functions and for severely ill-posed operators. On the other hand, the dense empirical minimizer is expected to work only for relatively smooth cases (the entropy integral has to

Received June 2008; revised June 2009.

[1]Supported by Deutsche Forschungsgemeinschaft under Project MA1026/8-1.

converge). However, because the minimization in the calculation of this estimator is not restricted to a $\delta$-net, we have available a larger toolbox of algorithms for finding (an approximation of) the minimizer of the empirical risk.

Let $(\mathbf{Y}, \mathcal{Y}, \nu)$ be a Borel space and let $A: L_2(\mathbf{R}^d) \to L_2(\mathbf{Y})$ be a linear operator, where $L_2(\mathbf{R}^d)$ is the space of square integrable functions $f: \mathbf{R}^d \to \mathbf{R}$ (with respect to the Lebesgue measure) and $L_2(\mathbf{Y})$ is the space of square integrable functions $g: \mathbf{Y} \to \mathbf{R}$ (with respect to measure $\nu$). In the density estimation model, we have i.i.d. observations

$$(1) \qquad Y_1, \ldots, Y_n \in \mathbf{Y}$$

with common density function $Af: \mathbf{Y} \to \mathbf{R}$, where $f: \mathbf{R}^d \to \mathbf{R}$ is a density function which we want to estimate. In the Gaussian white noise model, the observation is a realization of the process

$$(2) \qquad dY_n(y) = (Af)(y)\,dy + n^{-1/2}\,dW(y), \qquad y \in \mathbf{Y},$$

where $W(y)$ is the Brownian process on $\mathbf{Y}$, that is, for $h_1, h_2 \in L_2(\mathbf{Y})$, the random vector $(\int_{\mathbf{Y}} h_1\,dW, \int_{\mathbf{Y}} h_2\,dW)$ is a two-dimensional Gaussian random vector with zero mean, marginal variances $\|h_1\|_2^2, \|h_2\|_2^2$ and covariance $\int_{\mathbf{Y}} h_1 h_2\,d\nu$. (In our examples, $\mathbf{Y}$ is either the Euclidean space or the product of the real half-line with the unit sphere so that the existence of the Brownian process is guaranteed.) We want to estimate the signal function $f: \mathbf{R}^d \to \mathbf{R}$. The Gaussian white noise model is very useful for presenting the basic mathematical ideas in a transparent way. For the $\delta$-net estimator, the treatment is almost identical for the Gaussian white noise model and for the density estimation, but when we consider the dense empirical risk minimization, then, in the density estimation model, we need to use bracketing numbers and empirical entropies with bracketing, instead of the usual $L_2$-entropies. Our results for the Gaussian white noise model can also serve as a first step for obtaining analogous results for inverse problems in regression or in other statistical models.

The $L_2$-empirical risk is defined by

$$(3) \qquad \gamma_n(g) = \begin{cases} -2\int_{\mathbf{Y}} (Qg)\,dY_n + \|g\|_2^2, & \text{Gaussian white noise,} \\ -2n^{-1}\sum_{i=1}^{n} (Qg)(Y_i) + \|g\|_2^2, & \text{density estimation,} \end{cases}$$

where $Q$ is the adjoint of the inverse of $A$:

$$(4) \qquad \int_{\mathbf{R}^d} (A^{-1}h)g = \int_{\mathbf{Y}} h(Qg)\,d\nu$$

for $h \in L_2(\mathbf{Y})$, $g \in L_2(\mathbf{R}^d)$. The operator $Q = (A^{-1})^*$ has the domain $L_2(\mathbf{R}^d)$, similarly as $A$. Minimizing $\|\hat{f} - f\|_2^2$ with respect to estimators $\hat{f}$ is equivalent to minimizing $\|\hat{f} - f\|_2^2 - \|f\|_2^2$ and we have, in the Gaussian white noise

model,

$$\|\hat{f} - f\|_2^2 - \|f\|_2^2 = -2\int_{\mathbf{R}^d} f\hat{f} + \|\hat{f}\|_2^2$$

(5)
$$= -2\int_{\mathbf{Y}} (Af)(Q\hat{f})\,d\nu + \|\hat{f}\|_2^2$$

$$\approx -2\int_{\mathbf{Y}} (Q\hat{f})\,dY_n + \|\hat{f}\|_2^2 = \gamma_n(\hat{f}).$$

The usual least squares estimator is defined as a minimizer of the criterion

(6) $$\|A\hat{f} - Af\|_2^2 - \|Af\|_2^2 \approx -2\int_{\mathbf{Y}} (Ag)\,dY_n + \|Ag\|_2^2 \stackrel{\text{def}}{=} \tilde{\gamma}_n(g);$$

see, for example, O'Sullivan (1986). In density estimation, the log-likelihood empirical risk has been more common than the $L_2$-empirical risk and in the setting of inverse problems, the log-likelihood is defined as $\bar{\gamma}_n(g) = -n^{-1}\sum_{i=1}^n \log(Ag) \times (Y_i)$, analogously to (6). These alternative definitions of the empirical risk do not seem to lead to such an elegant theory as does the empirical risk in (3). The empirical risk in (3) has been used in deconvolution problems for projection estimators by Comte, Taupin and Rozenholc (2006).

We give upper bounds for the mean integrated squared error (MISE) of the estimators. The upper bounds characterize how the rates of convergence depend on the entropy of the underlying function class $\mathcal{F}$ and on smoothness properties of the operator $A$. Previously, such characterizations have been given (up to our knowledge) in inverse problems only for the case of estimating real-valued linear functionals $L$. In these cases, the rates of convergence are determined by the modulus of continuity of the functional $\omega(\epsilon) = \sup\{L(f) : f \in \mathcal{F}, \|Af\|_2 \leq \epsilon\}$; see Donoho and Low (1992). For the case of estimating the whole function with a global loss function, the rates of convergence depend on the size of the underlying function class in terms of the entropy and capacity; see Cencov (1972), Le Cam (1973), Ibragimov and Hasminskii (1980, 1981), Birgé (1983), Hasminskii and Ibragimov (1990), Yang and Barron (1999), Ibragimov (2004). $\delta$-net estimators were considered by, for example, van der Laan, Dudoit and van der Vaart (2004). These papers consider direct statistical problems. We show that for inverse statistical problems, the rate of convergence depends on the operator through the operator norm $\varrho(Q, \mathcal{F}_\delta)$ of $Q$, over a minimal $\delta$-net $\mathcal{F}_\delta$; see (8) for the definition of $\varrho(Q, \mathcal{F}_\delta)$. More precisely, the convergence rate $\psi_n$ of the $\delta$-net estimator is the solution to the equation

$$n\psi_n^2 = \varrho^2(Q, \mathcal{F}_{\psi_n})\log(\#\mathcal{F}_{\psi_n}),$$

where $\#\mathcal{F}_{\psi_n}$ is the cardinality of a minimal $\delta$-net. For direct problems, when $A$ is the identity operator, $\varrho(Q, \mathcal{F}_\delta) \asymp 1$. (We write $a_n \asymp b_n$ to mean that $0 < \liminf_{n\to\infty} a_n/b_n \leq \limsup_{n\to\infty} a_n/b_n < \infty$.) As examples of operators $A$, we consider the convolution operator and the Radon transform. For these operators, the estimators achieve the minimax rates of convergence over Sobolev classes.

The general framework for empirical risk minimization and the use of the empirical process machinery, including entropy bounds, for deriving optimal bounds seems to be new. Convolution and Radon transforms are discussed for illustrative purposes. These examples show that our results lead to optimal rates of convergence. As a new application, we introduce the estimation of additive models in inverse problems. A new type of oracle inequality is presented, which also gives the optimal rates of convergence in "anisotropic" inverse problems. For an extended version of this paper that also contains additional material, see Klemelä and Mammen (2009).

The paper is organized as follows. Section 2 gives an upper bound for the MISE of the $\delta$-net estimator. Section 3 gives a lower bound for the MISE of any estimator. Section 4 gives an upper bound for the MISE of the dense empirical risk minimizer. Section 5 proves that the $\delta$-net estimator achieves the optimal rate of convergence in the ellipsoidal framework and discusses this result for the case where $A$ is a convolution operator or the Radon transform. Furthermore, it contains an oracle inequality for additive models. Section 6 contains the proofs of the main results.

**2. $\delta$-net minimizer.** Let $\mathcal{F}$ be a set of densities or signal functions $f : \mathbf{R}^d \to \mathbf{R}$. Let $\mathcal{F}_\delta$ be a finite $\delta$-net of $\mathcal{F}$ in the $L_2$-metric, where $\delta > 0$. That is, for each $f \in \mathcal{F}$, there is a $\phi \in \mathcal{F}_\delta$ such that $\|f - \phi\|_2 \leq \delta$. Define the estimator $\hat{f}$ by

$$\hat{f} = \arg\min_{\phi \in \mathcal{F}_\delta} \gamma_n(\phi),$$

where $\gamma_n(\phi)$ is defined in (3). Typically, we would like to choose a $\delta$-net of minimal cardinality. We assume that $\mathcal{F}$ is bounded in the $L_2$-metric:

$$\sup_{g \in \mathcal{F}} \|g\|_2 \leq B_2, \tag{7}$$

where $0 < B_2 < \infty$.

Theorem 1 gives a bound for the mean integrated squared error of the estimator. We may identify the first term in the bound as a bias term and the second term as a variance term. The variance term depends on the operator norm of $Q$ over the $\delta$-net $\mathcal{F}_\delta$. We define this operator norm as

$$\varrho(Q, \mathcal{F}_\delta) = \max_{\phi, \phi' \in \mathcal{F}_\delta, \phi \neq \phi'} \frac{\|Q(\phi - \phi')\|_2}{\|\phi - \phi'\|_2}, \qquad \delta > 0, \tag{8}$$

where $Q$ is defined by (4). In the case of density estimation, we need the additional assumptions that $\varrho(Q, \mathcal{F}_\delta) \geq 1$ and that $A\mathcal{F}$ and $Q\mathcal{F}$ are bounded in the $L_\infty$ metric:

$$(9) \qquad \varrho(Q, \mathcal{F}_\delta) \geq 1, \qquad \sup_{f \in \mathcal{F}} \|Af\|_\infty \leq B_\infty, \qquad \sup_{f \in \mathcal{F}} \|Qf\|_\infty \leq B'_\infty,$$

where $0 < B_\infty, B'_\infty < \infty$.

THEOREM 1. *For the density estimation, we assume that (9) is satisfied. For $f \in \mathcal{F}$, we have that*

$$E\|\hat{f} - f\|_2^2 \leq C_1 \delta^2 + C_2 \frac{\varrho^2(Q, \mathcal{F}_\delta) \cdot (\log_e(\#\mathcal{F}_\delta) + 1)}{n},$$

*where*

$$(10) \qquad\qquad C_1 = (1 - 2\xi)^{-1}(1 + 2\xi),$$

$$(11) \qquad\qquad C_2 = (1 - 2\xi)^{-1} \xi C_\tau,$$

$$(12) \qquad\qquad C_\tau > 0$$

*and $\xi$ is such that*

$$(13) \qquad \begin{cases} C_\tau^{-1}(4B'_\infty/3 + \sqrt{2[8(B'_\infty)^2/9 + C_\tau B_\infty]}) \leq \xi < 1/2, \\ \qquad\qquad\qquad density\ estimation, \\ \sqrt{2/C_\tau} \leq \xi < 1/2, \qquad white\ noise. \end{cases}$$

A proof of Theorem 1 is given in Section 6.2.

REMARK 1. Theorem 1 shows that the $\delta$-net estimator achieves the rate of convergence $\psi_n$ when $\psi_n$ is the solution of the equation

$$(14) \qquad\qquad \psi_n^2 \asymp n^{-1} \varrho^2(Q, \mathcal{F}_{\psi_n}) \log(\#\mathcal{F}_{\psi_n}).$$

We calculate the rate under the assumptions that $\log(\#\mathcal{F}_\delta)$ and $\varrho(Q, \mathcal{F}_\delta)$ increase polynomially as $\delta$ decreases: we assume that one can find a $\delta$-net whose cardinality satisfies

$$\log(\#\mathcal{F}_\delta) = C\delta^{-b}$$

for some constants $b, C > 0$ and we assume that

$$\varrho(Q, \mathcal{F}_\delta) = C'\delta^{-a}$$

for some $a, C' > 0$ (in the direct case $a = 0$ and $C' = 1$). Then (14) can be written as $\psi_n^2 \asymp n^{-1} \psi_n^{-2a-b}$ and the rate of the $\delta$-net estimator is

$$(15) \qquad\qquad \psi_n \asymp n^{-1/[2(a+1)+b]}.$$

Let $\mathcal{F}$ be a set of $s$-smooth, $d$-dimensional functions such that $b = d/s$. Then the rate is $\psi_n \asymp n^{-s/[2(a+1)s+d]}$, which, for the direct case $a = 0$, gives the classical rate $\psi_n \asymp n^{-s/(2s+d)}$.

**3. A lower bound for MISE.** Theorem 2 gives a lower bound for the mean integrated squared error of any estimator when estimating densities or signal functions $f : \mathbf{R}^d \to \mathbf{R}$ in the function class $\mathcal{F}$. Theorem 2 also holds for nonlinear operators.

THEOREM 2. *Let $A$ be a possibly nonlinear operator. Assume that for each sufficiently small $\delta > 0$, we find a finite set $\mathcal{D}_\delta \subset \mathcal{F}$ for which*

$$\min\{\|f - g\|_2 : f, g \in \mathcal{D}_\delta, f \neq g\} \geq C_0 \delta \tag{16}$$

*and*

$$\begin{cases} \max\{\|f - g\|_2 : f, g \in \mathcal{D}_\delta\} \leq C_1 \delta, & \text{white noise,} \\ \max\{D_K(f, g) : f, g \in \mathcal{D}_\delta\} \leq C_1 \delta, & \text{density estimation,} \end{cases} \tag{17}$$

*where $D_K^2(f, g) = \int \log_e(f/g) f$ is the Kullback–Leibler distance and $C_0$, $C_1$ are positive constants. Let*

$$\varrho_K(A, \mathcal{D}_\delta) = \begin{cases} \dfrac{1}{\sqrt{2}} \max_{f, g \in \mathcal{D}_\delta, f \neq g} \dfrac{\|A(f - g)\|_2}{\|f - g\|_2}, & \text{white noise,} \\ \max_{f, g \in \mathcal{D}_\delta, f \neq g} \dfrac{D_K(Af, Ag)}{\|f - g\|_2}, & \text{density estimation.} \end{cases}$$

*Let $\psi_n$ be such that*

$$\log_e(\#\mathcal{D}_{\psi_n}) \succcurlyeq n \psi_n^2 \varrho_K^2(A, \mathcal{D}_{\psi_n}), \tag{18}$$

*where $a_n \succcurlyeq b_n$ means that $\liminf_{n \to \infty} a_n / b_n > 0$. Assume that*

$$\lim_{n \to \infty} n \psi_n^2 \varrho_K^2(A, \mathcal{D}_{\psi_n}) = \infty. \tag{19}$$

*Then*

$$\liminf_{n \to \infty} \psi_n^{-2} \inf_{\hat{f}} \sup_{f \in \mathcal{F}} E \|f - \hat{f}\|_2^2 > 0,$$

*where the infimum is taken over all estimators. That is, $\psi_n$ is a lower bound for the minimax rate of convergence.*

A proof of Theorem 2 is given in Section 6.3.

REMARK 2. Theorem 2 shows that one can get a lower bound $\psi_n$ for the rate of convergence by solving the equation

$$\psi_n^2 \varrho_K^2(A, \mathcal{D}_{\psi_n}) \asymp n^{-1} \log_e(\#\mathcal{D}_{\psi_n}). \tag{20}$$

The upper bound in Theorem 1 depends on the operator norm of $Q$, defined in (8), whereas the lower bound depends on the operator norm of $A$. Note, also, that the operator norm $\varrho(Q, \mathcal{F}_{\psi_n})$ is on the other side of the equation in (14) compared to the operator norm $\varrho_K(A, \mathcal{D}_{\psi_n})$ in the equation (20).

EMPIRICAL RISK MINIMIZATION 7

REMARK 3. In the density estimation case, one can easily check assumptions (17) and (19) if one assumes that the functions in $A\mathcal{D}_\delta$ are bounded and bounded away from 0. Then

$$(21) \qquad C' \cdot \|A(f-g)\|_2 \leq D_K(Af, Ag) \leq C \cdot \|A(f-g)\|_2$$

and (17) and (19) follow by the corresponding conditions with Hilbert norms instead of Kullback–Leibler distances.

**4. Dense minimizer.** The dense minimizer minimizes the empirical risk over the whole function class $\mathcal{F}$. In contrast to the $\delta$-net estimator, the minimization is not restricted to a $\delta$-net. We call this estimator a "dense minimizer" because it is defined as a minimizer over a possibly uncountable function class. The $\delta$-net estimator is more widely applicable: it may also be applied to estimate nonsmooth functions and it may be applied when the operator is severely ill-posed. The dense minimizer may only be applied for relatively smooth cases (the entropy integral has to converge). Because it works without a restriction to a $\delta$-net, we have available a larger toolbox of numerical algorithms that can be applied.

For a collection $\mathcal{F}$ of functions $f : \mathbf{R}^d \to \mathbf{R}$, the dense minimizer $\hat{f}$ is defined as a minimizer of the empirical risk over $\mathcal{F}$, up to $\epsilon_n > 0$:

$$\gamma_n(\hat{f}) \leq \inf_{g \in \mathcal{F}} \gamma_n(g) + \epsilon_n,$$

where $\gamma_n(\phi)$ is defined in (3). For clarity, we present separate theorems for the Gaussian white noise model and for the density estimation model. In both models, we make the assumption that the functions in $\mathcal{F}$ are bounded in the $L_2$-metric as in (7).

4.1. *Gaussian white noise.* Let $\mathcal{F}_\delta$, $\delta > 0$, be a $\delta$-net of $\mathcal{F}$, with respect to the $L_2$-norm. Define

$$(22) \qquad \varrho(Q, \mathcal{F}_\delta) = \max\left\{ \frac{\|Q(f-g)\|_2}{\|f-g\|_2} : f \in \mathcal{F}_\delta, g \in \mathcal{F}_{2\delta}, f \neq g \right\}, \qquad \delta > 0,$$

where $Q$ is the adjoint of the inverse of $A$, defined by (4). Define the entropy integral

$$(23) \qquad G(\delta) \stackrel{\text{def}}{=} \int_0^\delta \varrho(Q, \mathcal{F}_u) \sqrt{\log_e(\#\mathcal{F}_u)}\, du, \qquad \delta \in (0, B_2],$$

where $B_2$ is the $L_2$-bound defined by (7).

THEOREM 3. *Assume that:*

1. *the entropy integral in (23) converges;*

2. $G(\delta)/\delta^2$ is decreasing on the interval $(0, B_2]$;
3. $\varrho(Q, \mathcal{F}_\delta) = c\delta^{-a}$, where $0 \leq a < 1$ and $c > 0$;
4. $\lim_{\delta \to 0} G(\delta)\delta^{a-1} = \infty$;
5. $\delta \mapsto \varrho(Q, \mathcal{F}_\delta)\sqrt{\log_e(\#\mathcal{F}_\delta)}$ is decreasing on $(0, B_2]$.

Let $\psi_n$ be such that

(24) $$\psi_n^2 \geq Cn^{-1/2}G(\psi_n),$$

where $C$ is a positive constant, and assume that $\lim_{n \to \infty} n\psi_n^{2(1+a)} = \infty$. Then, for $f \in \mathcal{F}$,

$$E\|\hat{f} - f\|_2^2 \leq C'(\psi_n^2 + \epsilon_n)$$

for a positive constant $C'$, for sufficiently large $n$.

A proof of Theorem 3 is given in Section 6.4.

REMARK 4. Assumption 5 is a technical assumption which is used to replace a Riemann sum by an entropy integral. We prefer to write the assumptions in terms of the entropy integral in order to make them more readable.

REMARK 5. We may write $\varrho(Q, \mathcal{F}_\delta)$ in a simpler way when there exist minimal $\delta$-nets $\mathcal{F}_\delta$ which are nested: $\mathcal{F}_{2\delta} \subset \mathcal{F}_\delta$. We may then define, alternatively,

$$\varrho(Q, \mathcal{F}_\delta) = \max_{f,g \in \mathcal{F}_\delta, f \neq g} \frac{\|Q(f-g)\|_2}{\|f-g\|_2}.$$

REMARK 6. Theorems 3 and 4 show that the rate of convergence of the dense minimizer is the solution of the equation

(25) $$\psi_n^2 = n^{-1/2}G(\psi_n).$$

To get the optimal rate, the net $\mathcal{F}_\delta$ is chosen so that its cardinality is minimal. In the polynomial case, one can find a $\delta$-net whose cardinality satisfies

$$\log(\#\mathcal{F}_\delta) = C\delta^{-b}$$

for some constants $b, C > 0$ and the operator norm satisfies

$$\varrho(Q, \mathcal{F}_\delta) = C'\delta^{-a}$$

for some $a, C' > 0$. (In the direct case, $a = 0$ and $C' = 1$.) Thus, the entropy integral $G(\delta)$ is finite when $\int_0^\delta u^{-a-b/2}\,du < \infty$, which holds when

(26) $$a + b/2 < 1.$$

Then (25) leads to $\psi_n^2 \asymp n^{-1/2}\psi_n^{-a-b/2+1}$ and the rate of the dense minimization estimator is

$$\psi_n \asymp n^{-1/[2(a+1)+b]}. \tag{27}$$

This is the same rate as the rate of the $\delta$-net estimator given in (15). We have the following example. Let $\mathcal{F}$ be a set of $s$-smooth, $d$-dimensional functions such that $b = d/s$. Condition (26) may then be written as a condition for the smoothness index $s$: $s > d/[2(1-a)]$. When the problem is direct, then $a = 0$ and we have the classical condition $s > d/2$. The rate is $\psi_n \asymp n^{-s/[2(a+1)s+d]}$, which gives, for the direct case $a = 0$, the classical rate $\psi_n \asymp n^{-s/(2s+d)}$.

4.2. *Density estimation.* A $\delta$-bracketing net of $\mathcal{F}$ with respect to the $L_2$-norm is a set $\mathcal{F}_\delta = \{(g_j^L, g_j^U) : j = 1, \ldots, N_\delta\}$ of pairs of functions such that:

1. $\|g_j^L - g_j^U\|_2 \leq \delta$, $j = 1, \ldots, N_\delta$;
2. for each $g \in \mathcal{F}$, there exists $j = j(g) \in \{1, \ldots, N_\delta\}$ such that $g_j^L \leq g \leq g_j^U$.

Let us define $\mathcal{F}_\delta^L = \{g_j^L : j = 1, \ldots, N_\delta\}$ and $\mathcal{F}_\delta^U = \{g_j^U : j = 1, \ldots, N_\delta\}$. Further, define

$$\varrho_{\text{den}}(Q, \mathcal{F}_\delta) = \max\{\varrho(Q, \mathcal{F}_\delta^L, \mathcal{F}_\delta^U), \varrho(Q, \mathcal{F}_\delta^L, \mathcal{F}_{2\delta}^L)\}, \tag{28}$$

where

$$\varrho(Q, \mathcal{F}_\delta^L, \mathcal{F}_\delta^U) = \max\left\{\frac{\|Q(g^U - g^L)\|_2}{\|g^U - g^L\|_2} : g^L \in \mathcal{F}_\delta^L, g^U \in \mathcal{F}_\delta^U\right\}$$

and

$$\varrho(Q, \mathcal{F}_\delta^L, \mathcal{F}_{2\delta}^L) = \max\left\{\frac{\|Q(f - g)\|_2}{\|f - g\|_2} : f \in \mathcal{F}_\delta^L, g \in \mathcal{F}_{2\delta}^L, f \neq g\right\}$$

for $\delta > 0$. Define the entropy integral

$$G(\delta) \stackrel{\text{def}}{=} \int_0^\delta \varrho_{\text{den}}(Q, \mathcal{F}_u)\sqrt{\log_e(\#\mathcal{F}_u)}\, du, \qquad \delta \in (0, B_2], \tag{29}$$

where $B_2 = \sup_{f \in \mathcal{F}}\|f\|_2$.

THEOREM 4. *We make assumptions 1–5 of Theorem 3 [with operator norm $\varrho_{\text{den}}(Q, \mathcal{F}_\delta)$ in place of $\varrho(Q, \mathcal{F}_\delta)$] and, in addition, we assume that $\sup_{f \in \mathcal{F}}\|Af\|_\infty < \infty$, $\sup_{g \in \mathcal{F}_{B_2}^L \cup \mathcal{F}_{B_2}^U}\|Qg\|_\infty < \infty$ and that the operator $Q$ preserves positivity ($g \geq 0$ implies that $Qg \geq 0$). Let $\psi_n$ be such that*

$$\psi_n^2 \geq Cn^{-1/2}G(\psi_n) \tag{30}$$

for a positive constant $C$ and assume that $\lim_{n\to\infty} n\psi_n^{2(1+a)} = \infty$. Then, for $f \in \mathcal{F}$,
$$E\|\hat{f} - f\|_2^2 \le C'(\psi_n^2 + \epsilon_n)$$
for a positive constant $C'$, for sufficiently large $n$.

A proof of Theorem 4 is given in Section 6.5. An analogous discussion of optimal rates as in Remark 6 for the Gaussian white noise model also applies for dense density estimators.

**5. Examples of function spaces.** In Section 5.1, we consider ellipsoidal function spaces and in Section 5.2 we consider additive models and their generalizations.

5.1. *Ellipsoidal function spaces.* Since we are in the $L_2$-setting, it is natural to work in the sequence space; we define the function classes as ellipsoids. We shall apply singular value decompositions of the operators and wavelet-vaguelette systems in the calculation of the rates of convergence. In Section 5.1.1, we calculate the operator norms in the framework of singular value decompositions. In Section 5.1.2, we calculate the operator norms in the wavelet-vaguelettte framework. Section 5.1.3 derives the rate of convergence of the $\delta$-net estimator for the case of a convolution operator and the Radon transform, and the lower bound for the rate of convergence of any estimator.

5.1.1. *Singular value decomposition.* We assume that the underlying function space $\mathcal{F}$ consists of $d$-variate functions that are linear combinations of orthonormal basis functions $\phi_j$ with multi-index $j = (j_1, \ldots, j_d) \in \{0, 1, \ldots\}^d$. Define the ellipsoid and the corresponding collection of functions by

$$(31) \quad \Theta = \left\{\theta : \sum_{j_1=0,\ldots,j_d=0}^{\infty} a_j^2 \theta_j^2 \le L^2\right\}, \qquad \mathcal{F} = \left\{\sum_{j_1=0,\ldots,j_d=0}^{\infty} \theta_j \phi_j : \theta \in \Theta\right\}.$$

*$\delta$-net and $\delta$-packing set for polynomial ellipsoids.* We assume that there exist positive constants $C_1, C_2$ such that for all $j \in \{0, 1, \ldots\}^d$,

$$(32) \qquad C_1 \cdot |j|^s \le a_j \le C_2 \cdot |j|^s,$$

where $|j| = j_1 + \cdots + j_d$. In Klemelä and Mammen (2009), we construct a $\delta$-net $\Theta_\delta$ and a $\delta$-packing set $\Theta_\delta^*$ using the techniques of Kolmogorov and Tikhomirov (1961); see also Birman and Solomyak (1967). Since the construction is in the sequence space, we define the $\delta$-net and $\delta$-packing set of $\mathcal{F}$ by

$$(33) \quad \mathcal{F}_\delta = \left\{\sum_{j_1=0,\ldots,j_d=0}^{\infty} \theta_j \phi_j : \theta \in \Theta_\delta\right\}, \qquad \mathcal{D}_\delta = \left\{\sum_{j_1=0,\ldots,j_d=0}^{\infty} \theta_j \phi_j : \theta \in \Theta_\delta^*\right\}.$$

The set $\Theta_\delta$ is such that for $\theta \in \Theta_\delta$,

$$\theta_j = 0 \qquad \text{when } j \notin \{1, \ldots, M\}^d,$$

where

(34) $$M \asymp \delta^{-1/s}.$$

The set $\Theta_\delta^*$ is such that for all $\theta \in \Theta_\delta^*$,

(35) $$\theta_j = \theta_j^* \qquad \text{when } j \notin \{M^*, \ldots, M\}^d,$$

where $\theta^*$ is a fixed sequence with $\sum_{|j|\geq 0}^\infty a_j^2 \theta_j^{*2} = L^* < L$ and where $M^* = [M/2]$. Furthermore, it holds that

(36) $$\log(\#\Theta_\delta) \leq C\delta^{-d/s}, \qquad \log(\#\Theta_\delta^*) \geq C'\delta^{-d/s}.$$

*Operator norms.* We calculate the operator norms $\varrho(Q, \mathcal{F}_\delta)$ and $\varrho_K(A, \mathcal{D}_\delta)$ in the ellipsoidal framework, where $\mathcal{F}_\delta$ and $\mathcal{D}_\delta$ are defined in (33). We apply the singular value decomposition of $A$. We assume that the domain of $A$ is a separable Hilbert space $H$ with inner product $\langle \cdot, \cdot \rangle$. The underlying function space $\mathcal{F}$ satisfies $\mathcal{F} \subset H$. We denote by $A^*$ the adjoint of $A$. We assume that $A^*A$ is a compact operator on $H$ with eigenvalues $(b_j^2)$, $b_j > 0$, $j \in \{0, 1, \ldots\}^d$, with an orthonormal system of eigenfunctions $\phi_j$. We assume that there exist positive constants $q$ and $C_1, C_2$ such that for all $j \in \{0, 1, \ldots\}^d$,

(37) $$C_1 \cdot |j|^{-q} \leq b_j \leq C_2 \cdot |j|^{-q}.$$

Let $g, g'$ be elements of $\mathcal{F}_\delta$ or of $\mathcal{D}_\delta$, respectively. Write

$$g - g' = \sum_{j_1=1,\ldots,j_d=1}^\infty (\theta_j - \theta_j')\phi_j.$$

1. The functions $Q\phi_j$ are orthogonal and $\|Q\phi_j\|_2 = b_j^{-1}$. Indeed, $Q = (A^{-1})^*$ and thus

$$\langle Q\phi_j, Q\phi_l \rangle = \langle \phi_j, A^{-1}(A^{-1})^* \phi_l \rangle = b_l^{-2}\langle \phi_j, \phi_l \rangle,$$

where we have used the fact that

$$A^{-1}(A^{-1})^*\phi_l = A^{-1}(A^*)^{-1}\phi_l = (A^*A)^{-1}\phi_l = b_l^{-2}\phi_l.$$

Thus, for $g, g' \in \mathcal{F}_\delta$,

$$\|Q(g-g')\|_2^2 = \left\| \sum_{j_1=0,\ldots,j_d=0}^M (\theta_j - \theta_j')^2 Q\phi_j \right\|_2^2$$

$$\text{(38)} \qquad = \sum_{j_1=0,\ldots,j_d=0}^{M} (\theta_j - \theta'_j)^2 b_j^{-2}$$

$$\leq CM^{2q} \sum_{j_1=0,\ldots,j_d=0}^{M} (\theta_j - \theta'_j)^2,$$

where we have used (37) to infer that when $j \in \{0, \ldots, M\}^d$,
$$b_j^{-2} \leq C_1^{-2} \cdot |j|^{2q} \leq C_1^{-2} \cdot (dM)^{2q}.$$

On the other hand, $\|g - g'\|_2 = \sum_{j_1=0,\ldots,j_d=0}^{M} (\theta_j - \theta'_j)^2$. This gives the upper bound for the operator norm

$$\text{(39)} \qquad \varrho(Q, \mathcal{F}_\delta) \leq CM^q \leq C'\delta^{-q/s}$$

by the definition of $M$ in (34).

2. The functions $A\phi_j$ are orthogonal and $\|A\phi_j\|_2 = b_j$. Indeed,
$$\langle A\phi_j, A\phi_l \rangle = \langle \phi_j, A^*A\phi_l \rangle = b_l^2 \langle \phi_j, \phi_l \rangle.$$

Thus, for $g, g' \in \mathcal{D}_\delta$,

$$\text{(40)} \qquad \|A(g - g')\|_2^2 = \sum_{j_1=M^*,\ldots,j_d=M^*}^{M} (\theta_j - \theta'_j)^2 b_j^2.$$

This, together with calculations similar to those in (38), implies that

$$\text{(41)} \qquad C'\delta^{q/s} \leq \varrho_K(A, \mathcal{D}_\delta) \leq C\delta^{q/s}.$$

5.1.2. *Wavelet-vaguelette decomposition.* We assume that the underlying function space $\mathcal{F}$ consists of $d$-variate functions which are linear combinations of orthonormal wavelet functions $(\phi_{jk})$, where $j \in \{0, 1, \ldots\}$ and $k \in \{0, \ldots, 2^j - 1\}^d$. The $l_2$-body and the corresponding class of functions can now be defined as

$$\Theta = \left\{ \theta : \sum_j 2^{2sj} \sum_k |\theta_{jk}|^2 \leq L^2 \right\}, \qquad \mathcal{F} = \left\{ \sum_j \sum_k \theta_{jk}\phi_{jk} : \theta \in \Theta \right\},$$

where $s > 0$. We have already constructed a $\delta$-net and $\delta$-packing set for the $l_2$-bodies in (33). Now, this is done such that for $\theta \in \Theta_\delta$,
$$\theta_{jk} = 0 \qquad \text{when } j \geq J + 1,$$

where

$$\text{(42)} \qquad 2^J \asymp \delta^{-1/s}$$

and such that for $\theta \in \Theta^*_\delta$,
$$\theta_{jk} = \theta^*_{jk} \qquad \text{when } j \leq J^* \text{ or } j \geq J + 1,$$

where $\theta^*$ is a fixed sequence with $\sum_{j=0}^{\infty} \sum_k a_j^2 \theta^{*2}_{jk} = L^* < L$, and $J^* = J - 1$.

*Operator norms.* We can apply the wavelet-vaguelette decomposition, as defined in Donoho (1995), to calculate the operator norms $\varrho(Q, \mathcal{F}_\delta)$ and $\varrho_K(A, \mathcal{D}_\delta)$. We have available the following three sets of functions: $(\phi_{jk})_{jk}$ is an orthogonal wavelet basis and $(u_{jk})_{jk}$ and $(v_{jk})_{jk}$ are near-orthogonal sets:

$$\left\| \sum_{jk} a_{jk} u_{jk} \right\|_2 \asymp \|(a_{jk})\|_{l_2}, \qquad \left\| \sum_{jk} a_{jk} v_{jk} \right\|_2 \asymp \|(a_{jk})\|_{l_2},$$

where $a \asymp b$ means that there exist positive constants $C, C'$ such that $Cb \leq a \leq C'b$. The following quasi-singular relations hold:

$$A\phi_{jk} = \kappa_j v_{jk}, \qquad A^* u_{jk} = \kappa_j \phi_{jk},$$

where $\kappa_j$ are quasi-singular values. We assume that there exist positive constants $q$ and $C_1, C_2$ such that for all $j \in \{0, 1, \ldots\}$,

(43) $$C_1 \cdot 2^{-qj} \leq \kappa_j \leq C_2 \cdot 2^{-qj}.$$

1. Let $g, g' \in \mathcal{F}_\delta$. Write

$$g - g' = \sum_{j=0}^{J} \sum_k (\theta_{jk} - \theta'_{jk}) \phi_{jk}.$$

Since $Q = (A^{-1})^*$, it holds that $QA^* = (AA^{-1})^* = I$. Thus,

$$\langle Q\phi_{jk}, Q\phi_{j'k'} \rangle = \kappa_j^{-1} \kappa_{j'}^{-1} \langle QA^* u_{jk}, QA^* u_{j'k'} \rangle$$
$$= \kappa_j^{-1} \kappa_{j'}^{-1} \langle u_{jk}, u_{j'k'} \rangle.$$

This gives that

(44)
$$\|Q(g - g')\|_2^2 = \left\| \sum_{j=0}^{J} \sum_k (\theta_{jk} - \theta'_{jk}) Q\phi_{jk} \right\|_2^2$$
$$\asymp \sum_{j=0}^{J} \kappa_j^{-2} \sum_k (\theta_{jk} - \theta'_{jk})^2 \leq C 2^{2qJ} \sum_{j=0}^{J} \sum_k (\theta_{jk} - \theta'_{jk})^2,$$

where we have used (43) to infer that for $j \in \{0, \ldots, J\}$, it holds that
$$\kappa_j^{-2} \leq C_1^{-2} \cdot 2^{2qj} \leq C_1^{-2} \cdot 2^{2qJ}.$$

On the other hand, $\|g - g'\|_2^2 = \sum_{j=0}^{J} \sum_k (\theta_{jk} - \theta'_{jk})^2$. This gives the upper bound for the operator norm

$$\varrho(Q, \mathcal{F}_\delta) \leq C 2^{qJ} \leq C' \delta^{-q/s}$$

by the definition of $J$ in (42).

2. We have $\langle A\phi_{jk}, A\phi_{j'k'} \rangle = \kappa_j \kappa_{j'} \langle v_{jk}, v_{j'k'} \rangle$ and $(v_{jk})$ is a near-orthogonal set. Thus, similarly as in (44), we get

$$C' \delta^{q/s} \leq \varrho_K(A, \mathcal{D}_\delta) \leq C \delta^{q/s}.$$

5.1.3. *Rates of convergence.* We derive the rates of convergence for the $\delta$-net estimator when the operator is a convolution operator or the Radon transform. It is also shown that the lower bounds have the same order as the upper bounds. We will do this for the Gaussian white noise model.

*Convolution.* Let $A$ be a convolution operator: $Af = a * f$, where $a : \mathbf{R}^d \to \mathbf{R}$ is a known function and where $a * f(x) = \int_{\mathbf{R}^d} a(x-y) f(y) \, dy$ is the convolution of $a$ and $f$. For $j \in \{0, 1, \ldots\}^d$, $k \in K_j = \{k \in \{0,1\}^d : k_i = 0$, when $j_i = 0\}$, denote

$$\phi_{jk}(x) = \prod_{i=1}^{d} \sqrt{2}[(1-k_i)\cos(2\pi j_i x_i) + k_i \sin(2\pi j_i x_i)], \qquad x \in [0,1]^d.$$

The cardinality of $K_j$ is $2^{d-\alpha(j)}$, where $\alpha(j) = \#\{j_i : j_i = 0\}$. The collection $(\phi_{jk})$, $(j, k) \in \{0, 1, \ldots\}^d \times K_j$, is a basis for 1-periodic functions on $L_2([0,1]^d)$. When the convolution kernel $a$ is a 1-periodic function in $L_2([0,1]^d)$, then we can write

$$a(x) = \sum_{j_1=0,\ldots,j_d=0}^{\infty} \sum_{k \in K_j} b_{jk} \phi_{jk}(x).$$

The functions $\phi_{jk}$ are the singular functions of the operator $A$ and the values $b_{jk}$ are the corresponding singular values. We assume that the underlying function space is equal to

$$(45) \qquad \mathcal{F} = \left\{ \sum_{j_1=0,\ldots,j_d=0}^{\infty} \sum_{k \in K_j} \theta_{jk} \phi_{jk}(x) : (\theta_{jk}) \in \Theta \right\},$$

where

$$(46) \qquad \Theta = \left\{ \theta : \sum_{j_1=0,\ldots,j_d=0}^{\infty} \sum_{k \in K_j} a_{jk}^2 \theta_{jk}^2 \leq L^2 \right\}.$$

We give the rate of convergence of the $\delta$-net estimator and show that the estimator achieves the optimal rate of convergence. Optimal rates of convergence has been previously obtained for the convolution problem in various settings, in Ermakov (1989), Donoho and Low (1992), Koo (1993), Korostelev and Tsybakov (1993).

COROLLARY 1. *Let $\mathcal{F}$ be the function class as defined in (45). We assume that the coefficients of the ellipsoid (46) satisfy*

$$C_0 |j|^s \leq a_{jk} \leq C_1 |j|^s$$

for some $s > 0$ and $C_0, C_1 > 0$. We assume that the convolution filter $a$ is 1-periodic function in $L_2([0,1]^d)$ and that the Fourier coefficients of filter $a$ satisfy

$$C_2|j|^{-q} \leq b_{jk} \leq C_3|j|^{-q}$$

for some $q \geq 0$, $C_2, C_3 > 0$. Then

$$\limsup_{n \to \infty} n^{2s/(2s+2q+d)} \sup_{f \in \mathcal{F}} E_f \|\hat{f} - f\|_2^2 < \infty,$$

where $\hat{f}$ is the $\delta$-net estimator. Also,

$$\liminf_{n \to \infty} n^{2s/(2s+2q+d)} \inf_{\hat{g}} \sup_{f \in \mathcal{F}} E_f \|\hat{g} - f\|_2^2 > 0,$$

where the infimum is taken over any estimators $\hat{g}$.

PROOF. For the upper bound, we apply Theorem 1. Let $\mathcal{F}_\delta$ be the $\delta$-net of $\mathcal{F}$ as constructed in (33). We have shown in (39) that $\varrho(Q, \mathcal{F}_\delta) \leq C\delta^{-a}$, where $a = q/s$. We have stated in (36) that the cardinality of the $\delta$-net satisfies $\log(\#\mathcal{F}_\delta) \leq C\delta^{-b}$, where $b = d/s$. Thus, we may apply (15) to get the rate $\psi_n = n^{-1/(2(a+1)+b)} = n^{-s/(2s+2q+d)}$. This shows the upper bound.

For the lower bound, we apply Theorem 2. Assumption (16) holds because $\mathcal{D}_\delta$ in (33) is a $\delta$-packing set. Assumption (17) holds by the construction; see Klemelä and Mammen (2009). Assumptions (18) and (19) follow from (36) and (41). Thus the lower bound is proved. □

*Two-dimensional Radon transform.* We consider reconstructing a two-dimensional function $f$ from observations of its integrals over lines, that is, from its Radon transform. We suppose that $f \in L_1(D) \cap L_2(D)$, where $D = \{x \in \mathbf{R}^2 : \|x\| \leq 1\}$ is the unit disk in $\mathbf{R}^2$. We parametrize the lines by the length $u \in [0,1]$ of the perpendicular from the origin to the line and by the orientation $\phi \in [0, 2\pi)$ of this perpendicular. A common way to define the two-dimensional Radon transform is

$$(47) \quad Af(u, \phi) = \frac{\pi}{2\sqrt{1-u^2}} \int_{\sqrt{1-u^2}}^{\sqrt{1-u^2}} f(u\cos\phi - t\sin\phi, u\sin\phi + t\cos\phi)\,dt,$$

where $(u, \phi) \in \mathbf{Y} = [0, 1] \times [0, 2\pi]$. Now, the Radon transform is $\pi$ times the average of $f$ over the line segment that intersects $D$. We consider $Rf$ as the element of $L_2(\mathbf{Y}, \nu)$, where $\nu$ is the measure defined by $d\nu(u, \phi) = 2\pi^{-1}\sqrt{1-u^2}\,du\,d\phi$.

The singular value decomposition of the Radon transform can be found in Deans (1983). Let

$$\tilde{\phi}_{jk}(r, \theta) = \pi^{-1/2}(j+k+1)^{1/2} Z_{j+k}^{|j-k|}(r) e^{i(j-k)\theta},$$

$$(r, \theta) \in D = [0, 1] \times [0, 2\pi),$$

where $Z_a^b$ denotes the Zernike polynomial of degree $a$ and order $b$. Functions $\tilde{\phi}_{jk}$, $j,k=0,1,\ldots$, $(j,k)\neq(0,0)$, constitute an orthonormal complex-valued basis for $L_2(D)$. The corresponding orthonormal functions in $L_2(\mathbf{Y},\nu)$ are

$$\tilde{\psi}_{jk}(u,\phi) = \pi^{-1/2} U_{j+k}(u) e^{i(j-k)\phi}, \qquad (u,\phi) \in \mathbf{Y} = [0,1] \times [0,2\pi),$$

where $U_m(\cos\theta) = \sin((m+1)\theta)/\sin\theta$ are the Chebyshev polynomials of the second kind. We have

$$A\tilde{\phi}_{jk} = b_{jk}\tilde{\psi}_{jk},$$

where the singular values are

(48) $$b_{jk} = \pi^{-1}(j+k+1)^{-1/2}.$$

The complex basis identifies the equivalent real orthonormal basis as follows:

$$\phi_{jk} = \begin{cases} \sqrt{2}\,\mathrm{Re}(\tilde{\phi}_{jk}), & \text{if } j > k, \\ \tilde{\phi}_{jk}, & \text{if } j = k, \\ \sqrt{2}\,\mathrm{Im}(\tilde{\phi}_{jk}), & \text{if } j < k. \end{cases}$$

We assume that the underlying function space is equal to

(49) $$\mathcal{F} = \left\{ \sum_{j_1=0,j_2=0,(j_1,j_2)\neq(0,0)}^{\infty} \theta_{j_1 j_2} \phi_{j_1 j_2}(x) : (\theta_{j_1 j_2}) \in \Theta \right\},$$

where

(50) $$\Theta = \left\{ \theta : \sum_{j_1=0,j_2=0,(j_1,j_2)\neq(0,0)}^{\infty} a_{j_1 j_2}^2 \theta_{j_1 j_2}^2 \leq L^2 \right\}.$$

We give the rate of convergence of the $\delta$-net estimator and show that the estimator achieves the optimal rate of convergence. Optimal rates of convergence have previously been obtained in Johnstone and Silverman (1990), Korostelev and Tsybakov (1991), Donoho and Low (1992), Korostelev and Tsybakov (1993).

COROLLARY 2. *Let $\mathcal{F}$ be the function class as defined in (49). We assume that the coefficients of the ellipsoid (50) satisfy*

$$C_0 |j|^s \leq a_{jk} \leq C_1 |j|^s$$

*for some $s > 0$ and $C_0, C_1 > 0$. Then, for $d = 2$,*

$$\limsup_{n\to\infty} n^{2s/(2s+2d-1)} \sup_{f\in\mathcal{F}} E_f \|\hat{f} - f\|_2^2 < \infty,$$

*where $\hat{f}$ is the $\delta$-net estimator. Also,*

$$\liminf_{n\to\infty} n^{2s/(2s+2d-1)} \inf_{\hat{g}} \sup_{f\in\mathcal{F}} E_f \|\hat{g} - f\|_2^2 > 0,$$

*where the infimum is taken over any estimator $\hat{g}$.*

PROOF. For the upper bound, we apply Theorem 1. Let $\mathcal{F}_\delta$ be the $\delta$-net of $\mathcal{F}$ as constructed in (33). We have shown in (39) that

$$\varrho(Q, \mathcal{F}_\delta) \leq C\delta^{-a},$$

where $a = q/s$ and $q = 1/2$ [so that $a = (d-1)/(2s)$] since the singular values are given in (48). We have stated in (36) that the cardinality of the $\delta$-net satisfies

$$\log(\#\mathcal{F}_\delta) \leq C\delta^{-b},$$

where $b = d/s$. Thus, we can apply (15) to get the rate

$$\psi_n = n^{-s/(2s+2d-1)}.$$

The upper bound is proved. For the lower bound, we apply Theorem 2 similarly as in the proof of Corollary 1. $\square$

5.2. *Additive models.* In this section, we will show that our approach can be used to prove oracle results for additive models. In additive models, the unknown function $f : \mathbf{R}^d \to \mathbf{R}$ is assumed to have an additive decomposition $f(x) = f_1(x_1) + \cdots + f_d(x_d)$ with unknown additive components $f_j : \mathbf{R} \to \mathbf{R}$, $j = 1, \ldots, d$. We compare this model with theoretical oracle models where only one component function $f_r$ is unknown, the other functions $f_j$ ($j \neq r$) being known. We will show below that the function $f$ can be estimated with the same rate of convergence as in the oracle model that has the slowest rate of convergence. In particular, if the rate of convergence is the same in all oracle models, then the rate in the additive model remains the same. This is a well-known fact for classical additive regression models; see, for example, Stone (1985). It efficiently avoids the curse of dimensionality, in contrast to the full-dimensional nonparametric model. Furthermore, it is practically important because it allows a flexible and nicely interpretable model for regression with high-dimensional covariates; see, for example, Hastie and Tibshirani (1990) for a discussion of the additive and related models. Thus, our result will generalize the oracle result for additive models of Stone (1985) to inverse problems. For a theoretical discussion, we will first use a slightly more general framework. We will later return to additive models.

5.2.1. *Abstract setting.* We assume that the function class $\mathcal{F}$ is a subset of the direct sum of spaces $\mathcal{F}_1, \ldots, \mathcal{F}_p$. All spaces contain functions from $\mathbf{R}^d \to \mathbf{R}$. At this stage, we do not assume that functions in $\mathcal{F}_j$ ($j = 1, \ldots, p$) depend only on the argument $x_j$. Examples of this more general set-up are sums of smooth functions and indicator functions of convex sets or of sets with smooth boundary. We assume that a finite $\delta$-net $\mathcal{F}_\delta$ of $\mathcal{F}$ is a subset of the direct sum $\mathcal{F}_{1,\delta} \oplus \cdots \oplus \mathcal{F}_{p,\delta}$, where $\mathcal{F}_{j,\delta}$ are finite subsets of $\mathcal{F}_j$. We denote the number of elements of $\mathcal{F}_{j,\delta}$ by $\exp(\lambda_j)$. Furthermore, we write $\rho_j = \rho(Q, \mathcal{F}_{j,\delta})$. We make the following, essential, geometrical assumption:

$$\|f_1 + \cdots + f_p\|_2^2 \geq c \sum_{j=1}^p \|f_j\|_2^2 \tag{51}$$

for a positive constant $c > 0$. For the $\delta$-net minimizer $\hat{f}$ over the $\delta$-net $\mathcal{F}_\delta$, we get the following result in the white noise model. (An additive model for density estimation would not make much sense.)

THEOREM 5. *We make assumption (51). In the white noise model, the following bound holds for the $\delta$-net minimizer $\hat{f}$, for $f \in \mathcal{F}$,*

$$E(\|\hat{f} - f\|_2^2) \leq 3\delta^2 + 32c^{-1}n^{-1}\left[\sum_{j=1}^p \rho_j^2 \lambda_j + \left(\sum_{j=1}^p \rho_j\right)^2\right].$$

A proof of Theorem 5 is given in Section 6.6.

5.2.2. *Application to additive models.* We now apply Theorem 5 in order to discuss additive models $f(x) = f_1(x_1) + \cdots + f_d(x_d)$. In $L_2(\mathbf{R}^d)$, we have $\|f_1 + \cdots + f_d\|_2^2 = \sum_{j=1}^d \|f_j\|_2^2$ if the functions $f_j$ are normed such that $\int f_j(x_j)\,dx_j = 0$. Thus, (51) holds trivially. Assumption (51) also holds in other $L_2$-spaces with dominating measure differing from the Lebesgue measure. A discussion of condition (51) for these classes can be found in, for example, Mammen, Linton and Nielsen (1999); also, see Bickel et al. (1993). Such $L_2$-spaces naturally arise in additive regression models. For a white noise model, they arise if one assumes an additive model for transformed covariables. We assume that for the models $\mathcal{F}_j$, one can find $\delta_j$-nets $\mathcal{F}_{j,\delta_j}$ such that choosing $\delta_j = \psi_{n,j}$ with

$$\psi_{n,j}^2 \asymp n^{-1}\rho^2(Q, \mathcal{F}_{j,\psi_{n,j}})\log(\#\mathcal{F}_{j,\psi_{n,j}})$$

gives a rate-optimal $\delta$-net minimizer in the model $\mathcal{F}_j$. Now, $\mathcal{F}_\delta = \mathcal{F}_{1,\delta_1} \oplus \cdots \oplus \mathcal{F}_{d,\delta_d}$ is a $\delta$-net of $\mathcal{F}$ with $\delta = \sum_{j=1}^d \delta_j$. From Theorem 5, we get that the $\delta$-net minimizer $\hat{f}$ over the net $\mathcal{F}_\delta$ achieves the rate $O(\psi_n)$ with $\psi_n = \max_{1 \leq j \leq d} \psi_{n,j}$. This is just the type of result we called an oracle result at the beginning of this section.

In general, the oracle result does not follow from Theorem 1. The application of Theorem 1 leads to an assumption of the type

$$n^{-1} \max_{1 \leq j \leq d} \rho^2(Q, \mathcal{F}_{j,\psi_{n,j}}) \times \max_{1 \leq j \leq d} \log(\#\mathcal{F}_{j,\psi_{n,j}}) = O(\psi_n^2),$$

whereas Theorem 5 only requires that
$$n^{-1} \max_{1 \le j \le d} [\rho^2(Q, \mathcal{F}_{j,\psi_{n,j}}) \log(\#\mathcal{F}_{j,\psi_{n,j}})] = O(\psi_n^2).$$

This can make a big difference. First, the entropy numbers of the additive classes $\mathcal{F}_j$ may differ. Furthermore, the operator $Q$ may act quite differently on the spaces $\mathcal{F}_j$.

5.2.3. *Ellipsoidal spaces and convolution.* As an example, we now assume that the underlying function space is $\mathcal{F} = \mathcal{F}_1 \oplus \cdots \oplus \mathcal{F}_d$, where
$$\mathcal{F}_k = \left\{ \sum_{j=0}^{\infty} \theta_{kj} \phi_{kj} : \theta_{k\cdot} \in \Theta_{s_k, L_k} \right\}$$

for basis functions $\phi_{kj} : [0,1] \to \mathbf{R}$ and the ellipsoids are defined by

(52) $$\Theta_{s_k, L_k} = \left\{ \theta_{k\cdot} : \sum_{j=0}^{\infty} a_{kj}^2 \theta_{kj}^2 \le L_k^2 \right\}, \qquad k = 1, \ldots, d,$$

where we assume that there exist positive constants $C_1, C_2$ such that for all $j \in \{0, 1, \ldots\}$,

(53) $$C_1 \cdot j^{s_k} \le a_{kj} \le C_2 \cdot j^{s_k}.$$

Let $A$ be a convolution operator: $Af = a * f$, where $a : \mathbf{R}^d \to \mathbf{R}$ is a known function. Then
$$Af = A_1 f_1 + \cdots + A_d f_d,$$
where $f(x) = f_1(x_1) + \cdots + f_d(x_d)$ and
$$A_k f_k(x_k) = \int_{[0,1]^d} f_k(x_k - y_k) a_k(y_k) \, dy_k,$$
where $a_k(y_k) = \int_{[0,1]^d} a(y) \prod_{l=1, l \ne k}^{d} dy_l$ is the $k$th marginal function of $a$. We can decompose $Q$ accordingly:
$$Qg = Q_1 g_1 + \cdots + Q_d g_d.$$
Operators $A_j$ and $Q_j$ are restrictions of $A$ and $Q$ to $\mathcal{F}_j$. We apply the singular value decomposition for $A_k$. Let
$$\phi_{kj}(t) = \sqrt{2} \cos(2\pi j t), \qquad t \in [0,1],$$
where $j = 1, 2, \ldots$ and $\phi_0(t) = I_{[0,1]}(t)$. The collection $(\phi_{kj})$, $j = 0, 1, \ldots$, is a basis for 1-periodic functions on $L_2([0,1])$. When $a_k$ are 1-periodic functions in $L_2([0,1])$, we can write
$$a_k(x_k) = \sum_{j=0}^{\infty} b_{kj} \phi_{kj}(x_k).$$

The functions $\phi_{kj}$ are the singular functions of the operator $A_k$ and the values $b_{kj}$ are the corresponding singular values. We give the rate of convergence of the $\delta$-net estimator and show that the estimator achieves the optimal rate of convergence.

COROLLARY 3. *Let $\mathcal{F} = \mathcal{F}_1 \oplus \cdots \oplus \mathcal{F}_d$. We assume that the coefficients of the ellipsoid satisfy (53). We assume that $a_k$ are 1-periodic functions in $L_2([0,1])$ and that the Fourier coefficients of $a_k$ satisfy*

$$C_2 j^{-q_k} \leq b_{kj} \leq C_3 j^{-q_k}$$

*for some $q_k \geq 0$, $C_2, C_3 > 0$. Then, in the white noise model,*

$$\limsup_{n \to \infty} n^a \sup_{f \in \mathcal{F}} E_f \|\hat{f} - f\|_2^2 < \infty,$$

*where $\hat{f}$ is the $\delta$-net estimator and*

$$a = \min_{k=1,\ldots,d} \frac{2s_k}{2s_k + 2q_k + 1}.$$

*Also,*

$$\liminf_{n \to \infty} n^a \inf_{\hat{g}} \sup_{f \in \mathcal{F}} E_f \|\hat{g} - f\|_2^2 > 0,$$

*where the infimum is taken over any estimators $\hat{g}$ in the white noise model.*

PROOF. For the upper bound, we apply Theorem 5. As in Section 5.1.1, we can find $\delta$-nets $\mathcal{F}_{k,\delta}$ for $\mathcal{F}_k$ whose cardinality is bounded by $\log(\#\mathcal{F}_{k,\delta}) \leq C\delta^{-1/s_k}$ and $\varrho(Q_k, \mathcal{F}_{k,\delta}) \leq C\delta^{-q_k/s_k}$. The upper bound of Theorem 5 gives as the rate the maximum of the component rates $n^{-2s_k/(2s_k+2q_k+1)}$. For the lower bound, we apply the lower bound of Corollary 1 in the case $d=1$ and the fact that one cannot do better in the additive model than in the model that has only one component. □

## 6. Proofs.

6.1. *A preliminary lemma.* We prove that the theoretical error of a minimization estimator may be bounded by the optimal theoretical error and an additional stochastic term.

LEMMA 1. *Let $\mathcal{C} \subset L_2(\mathbf{R}^d)$. Let $\hat{f} \in \mathcal{C}$ be such that*

(54) $$\gamma_n(\hat{f}) \leq \inf_{g \in \mathcal{C}} \gamma_n(g) + \varepsilon,$$

*where $\varepsilon \geq 0$. Then, for each $f^0 \in \mathcal{C}$,*

$$\|\hat{f} - f\|_2^2 \leq \|f^0 - f\|_2^2 + \varepsilon + 2\nu_n[Q(\hat{f} - f^0)],$$

where $f$ is the true density or the true signal function and $\nu_n(g)$ is the centered empirical operator:

$$(55) \quad \nu_n(g) = \begin{cases} \int g\, dY_n - \int_{\mathbf{Y}} g(Af), & \text{white noise model,} \\ n^{-1} \sum_{i=1}^n g(Y_i) - \int_{\mathbf{Y}} g(Af), & \text{density estimation,} \end{cases}$$

where $g: \mathbf{R}^d \to \mathbf{R}$.

PROOF. We have, for $g = \hat{f}$, $g = f^0$,

$$\|g - f\|_2^2 - \gamma_n(g)$$
$$= \begin{cases} \|f\|_2^2 - 2\int_{\mathbf{R}^d} fg + 2\int (Qg)\, dY_n, & \text{white noise model,} \\ \|f\|_2^2 - 2\int_{\mathbf{R}^d} fg + 2n^{-1}\sum_{i=1}^n (Qg)(Y_i), & \text{density estimation.} \end{cases}$$

We have $\int_{\mathbf{R}^d} fg = \int_{\mathbf{Y}} (Af)(Qg)$. Thus,

$$(56) \quad \|\hat{f} - f\|_2^2 - \gamma_n(\hat{f}) + \gamma_n(f^0) - \|f^0 - f\|_2^2 = 2\nu_n[Q(\hat{f} - f^0)].$$

Thus,

$$\|\hat{f} - f\|_2^2 - \|f^0 - f\|_2^2 = \|\hat{f} - f\|_2^2 - \gamma_n(\hat{f}) + \gamma_n(\hat{f}) - \|f^0 - f\|_2^2$$
$$(57) \qquad \qquad \leq \|\hat{f} - f\|_2^2 - \gamma_n(\hat{f}) + \gamma_n(f^0) + \varepsilon - \|f^0 - f\|_2^2$$
$$(58) \qquad \qquad = 2\nu_n[Q(\hat{f} - f^0)] + \varepsilon.$$

In (57), we applied (54) and in (58), we applied (56). □

6.2. *Proof of Theorem 1.* Let $f \in \mathcal{F}$ be the true density. Let $\phi^0 \in \mathcal{F}_\delta$. Define $\zeta = C_1 \|\phi^0 - f\|_2^2 + C_2 n^{-1} \varrho^2(Q, \mathcal{F}_\delta) \log_e(\#\mathcal{F}_\delta)$, where $C_1$ is defined in (10) and $C_2$ is defined in (11). We have that

$$E\|\hat{f} - f\|_2^2$$
$$(59) \quad = \int_0^\infty P(\|\hat{f} - f\|_2^2 > t)\, dt \leq \zeta + \int_\zeta^\infty P(\|\hat{f} - f\|_2^2 > t)\, dt$$
$$= \zeta + C_2 n^{-1} \varrho^2(Q, \mathcal{F}_\delta) \int_0^\infty P(\|\hat{f} - f\|_2^2 > C_2 n^{-1} \varrho^2(Q, \mathcal{F}_\delta) t + \zeta)\, dt.$$

Define $\tau_n = C_\tau n^{-1} \varrho^2(Q, \mathcal{F}_\delta)(\log_e(\#\mathcal{F}_\delta) + t)$, where $C_\tau$ is defined in (12). Then

$$P(\|\hat{f} - f\|_2^2 > C_2 n^{-1} \varrho^2(Q, \mathcal{F}_\delta) t + \zeta)$$

$$= P(\|\hat{f} - f\|_2^2 > C_1\|\phi^0 - f\|_2^2 + C_2 C_\tau^{-1}\tau_n)$$

(60)
$$= P((1 - 2\xi)^{-1}\|\hat{f} - f\|_2^2$$
$$> 2\xi(1 - 2\xi)^{-1}\|\hat{f} - f\|_2^2 + C_1\|\phi^0 - f\|_2^2 + C_2 C_\tau^{-1}\tau_n)$$
$$= P(\|\hat{f} - f\|_2^2 > 2\xi\|\hat{f} - f\|_2^2 + (1 + 2\xi)\|\phi^0 - f\|_2^2 + \xi\tau_n).$$

We have, by Lemma 1, that $\|\hat{f} - f\|_2^2 \leq \|\phi^0 - f\|_2^2 + 2\nu_n[Q(\hat{f} - \phi^0)]$. This implies that

$$P(\|\hat{f} - f\|_2^2 > C_2 n^{-1}\varrho^2(Q, \mathcal{F}_\delta)t + \zeta)$$
$$= P(\nu_n[Q(\hat{f} - \phi^0)] > \xi\|\hat{f} - f\|_2^2 + \xi\|\phi^0 - f\|_2^2 + \xi\tau_n/2)$$
(61)
$$= P(\nu_n[Q(\hat{f} - \phi^0)] > w(\hat{f})\xi)$$
$$\leq P\left(\max_{\phi \in \mathcal{F}_\delta, \phi \neq \phi^0} \frac{\nu_n[Q(\phi - \phi^0)]}{w(\phi)} > \xi\right)$$
$$\stackrel{\text{def}}{=} P_{\max},$$

where $w(\phi) = \|\phi - f\|_2^2 + \|\phi^0 - f\|_2^2 + \tau_n/2$. We will prove that

(62) $$P_{\max} \leq \exp(-t).$$

Together with (59) and (61), this proves the theorem.

*Proof of (62).* Define

$$\mathcal{G} = \left\{\frac{Q(\phi - \phi^0)}{w(\phi)} : \phi \in \mathcal{F}_\delta, \phi \neq \phi^0\right\}.$$

We have that

(63) $$P_{\max} \leq \sum_{g \in \mathcal{G}} P(\nu_n(g) > \xi).$$

Also, $w(\phi) \geq \frac{1}{2}(\|\phi - \phi^0\|_2^2 + \tau_n) \geq \|\phi - \phi^0\|_2 \tau_n^{1/2}$ and thus

(64) $$v_0 \stackrel{\text{def}}{=} \max_{g \in \mathcal{G}} \|g\|_2^2 \leq \frac{1}{\tau_n} \max_{\phi \in \mathcal{F}_\delta, \phi \neq \phi^0} \frac{\|Q(\phi - \phi_0)\|_2^2}{\|\phi - \phi_0\|_2^2} = \frac{\varrho^2(Q, \mathcal{F}_\delta)}{\tau_n}.$$

*Gaussian white noise.* When $W \sim N(0, \sigma^2)$, we have $P(W > \xi) \leq 2^{-1} \times \exp\{-\xi^2/(2\sigma^2)\}$ for $\xi > 0$; see, for example, Dudley (1999), Proposition 2.2.1. We have that $\nu_n(g) \sim N(0, n^{-1}\|g\|_2^2)$. Thus,

$$P(\nu_n(g) > \xi) \leq 2^{-1}\exp\left\{-\frac{n\xi^2}{2v_0}\right\} \leq 2^{-1}\exp\left\{-\frac{n\tau_n\xi^2}{2\varrho^2(Q, \mathcal{F}_\delta)}\right\}.$$

Defining $C_\xi \stackrel{\text{def}}{=} \xi^2 C_\tau/2$, we get that

$$P_{\max} \leq \#\mathcal{F}_\delta \cdot \exp\left\{-\frac{n\tau_n\xi^2}{2\varrho^2(Q,\mathcal{F}_\delta)}\right\} = \#\mathcal{F}_\delta \cdot \exp\{-C_\xi[\log_e(\#\mathcal{F}_\delta) + t]\}$$
$$\leq \exp(-t)$$

since $C_\xi \geq 1$, by the choice of $\xi$.

*Density estimation.* Define $v = \sup_{g\in\mathcal{G}} \text{Var}_f(g(Y_1))$ and $b = \sup_{g\in\mathcal{G}} \|g\|_\infty$. We have that

$$(65) \qquad v \leq \|Af\|_\infty v_0 \leq B_\infty \frac{\varrho^2(Q,\mathcal{F}_\delta)}{\tau_n}$$

by (64). Also, $w(\phi) \geq \tau_n/2$ and thus, because of $\varrho(Q,\mathcal{F}_\delta) \geq 1$, we have that

$$(66) \qquad b \leq 2B'_\infty \frac{2}{\tau_n} \leq 4B'_\infty \frac{\varrho^2(Q,\mathcal{F}_\delta)}{\tau_n}.$$

By applying Bernstein's inequality, we get, with (65) and (66), that

$$P(\nu_n(g) > \xi) \leq \exp\left\{\frac{-n\xi^2}{2(v + \xi b/3)}\right\}$$
$$\leq \exp\left\{\frac{-n\xi^2 \tau_n}{2\varrho^2(Q,\mathcal{F}_\delta)(B_\infty + 4B'_\infty\xi/3)}\right\}.$$

Continuing from (63),

$$P_{\max} \leq \#\mathcal{F}_\delta \cdot \exp\left\{\frac{-n\xi^2 \tau_n}{2\varrho^2(Q,\mathcal{F}_\delta)(B_\infty + 4B'_\infty\xi/3)}\right\}$$
$$= \#\mathcal{F}_\delta \cdot \exp\{-C_\xi[\log_e(\#\mathcal{F}_\delta) + t]\} \leq \exp(-t),$$

where

$$C_\xi \stackrel{\text{def}}{=} \frac{\xi^2 C_\tau}{2(B_\infty + 4B'_\infty\xi/3)}$$

and $C_\xi \geq 1$, by the choice of $\xi$. We have proven (62) and thus the theorem.

6.3. *Proof of Theorem 2.* To prove Theorem 2, we follow the approach of Hasminskii and Ibragimov (1990) and use Theorem 6 in Tsybakov (1998), which gives us the following lemma. For a proof, see Klemelä and Mammen (2009).

LEMMA 2. *Let $\mathcal{D} \subset \mathcal{F}$ be a finite set for which*

$$(67) \qquad \min\{\|f - g\|_2 : f, g \in \mathcal{D}, f \neq g\} \geq \delta,$$

where $\delta > 0$. Assume that for some $f_0 \in \mathcal{D}$ and for all $f \in \mathcal{D} \setminus \{f_0\}$,

$$P_{Af}^{(n)}\left(\frac{dP_{Af_0}^{(n)}}{dP_{Af}^{(n)}} \leq \tau\right) \leq \alpha, \tag{68}$$

for some $0 < \alpha < 1$, $\tau > 0$. Here, $P_{Af}^{(n)}$ is the product measure corresponding to the density $Af$ in the density estimation model, and in the Gaussian white noise model, $P_{Af}^{(n)}$ is the measure of the process $Y_n$ in (2). It then holds that

$$\inf_{\hat{f}} \sup_{f \in \mathcal{F}} E_{Af}\|f - \hat{f}\|_2^2 \geq \frac{\delta^2}{4}(1-\alpha)\frac{\tau(N_\delta - 1)}{1 + \tau(N_\delta - 1)},$$

where $N_\delta = \#\mathcal{D} \geq 2$. Here, the infimum is taken over all estimators (either in the density estimation model or in the Gaussian white noise model).

PROOF OF THEOREM 2. For $f, f_0 \in \mathcal{D}_{\psi_n}$, $f \neq f_0$,

$$P_{Af}^{(n)}\left(\frac{dP_{Af_0}^{(n)}}{dP_{Af}^{(n)}} \leq \tau\right)$$

$$\leq (\log \tau^{-1})^{-1} D_K^2(P_{Af}^{(n)}, P_{Af_0}^{(n)}) \tag{69}$$

$$= \begin{cases} (\log \tau^{-1})^{-1} n D_K^2(Af, Af_0), & \text{density estimation,} \\ (\log \tau^{-1})^{-1} \frac{n}{2}\|Af - Af_0\|_2^2, & \text{Gaussian white noise,} \end{cases} \tag{70}$$

where, in (69), we applied Markov's inequality and for the Gaussian white noise model, in (70), we applied the fact that under $P_{Af}^{(n)}$,

$$\frac{dP_{Af}^{(n)}}{dP_{Af_0}^{(n)}} = \exp\{n^{1/2}\sigma Z + n\sigma^2/2\},$$

where $Z \sim N(0, 1)$ and $\sigma = \|Af - Af_0\|_2$. When we choose

$$\tau = \tau_n = \exp\{-\alpha^{-1}n[C_1 \varrho_K(A, \mathcal{D}_{\psi_n})\psi_n]^2\}$$

for $0 < \alpha < 1$, we get, by applying assumption (17), that

$$P_{Af}^{(n)}\left(\frac{dP_{Af_0}^{(n)}}{dP_{Af}^{(n)}} \leq \tau\right) \leq (\log \tau^{-1})^{-1} n \varrho_K^2(A, \mathcal{D}_{\psi_n})\|f - f_0\|_2^2$$

$$\leq (\log \tau^{-1})^{-1} n[\varrho_K(A, \mathcal{D}_{\psi_n}) C_1 \psi_n]^2 = \alpha. \tag{71}$$

By applying Lemma 2, assumption (16) and (71), we get the lower bound

$$\inf_{\hat{f}} \sup_{f \in \mathcal{D}_{\psi_n}} \|f - \hat{f}\|_2^2 \geq \frac{(C_0\psi_n)^2}{4}(1-\alpha)\frac{\tau_n(N_{\psi_n} - 1)}{1 + \tau_n(N_{\psi_n} - 1)}, \tag{72}$$

where $N_{\psi_n} = \#\mathcal{D}_{\psi_n}$. Let $n$ be so large that $\log_e N_{\psi_n} \geq C_2^2 n \varrho_K^2(A, \mathcal{D}_{\psi_n}) \psi_n^2$, where $C_2 > C_1$. This is possible by (18). Then

$$\tau_n N_{\psi_n} = \exp\{\log_e N_{\psi_n} - \alpha^{-1} n [C_1 \varrho_K(A, \mathcal{D}_{\psi_n}) \psi_n]^2\}$$
$$\geq \exp\{n \varrho_K^2(A, \mathcal{D}_{\psi_n}) \psi_n^2 [C_2^2 - \alpha^{-1} C_1^2]\} \to \infty$$

as $n \to \infty$, where we apply (19) and choose $\alpha$ so that $C_2^2 - \alpha^{-1} C_1^2 > 0$, that is, $(C_1/C_2)^2 < \alpha < 1$. Then $\lim_{n \to \infty} \tau_n(N_{\psi_n} - 1)/[1 + \tau_n(N_{\psi_n} - 1)] = 1$ and the theorem follows from (72). $\square$

6.4. *Proof of Theorem 3.* Let $\zeta = C_1 \epsilon_n + C_2 \psi_n^2$, where $C_1 = (1 - 2\xi)^{-1}$, $C_2 = 1 - 2\xi$, $0 < \xi \leq (3 - \sqrt{5})/4$. We have that

$$(73) \quad E\|\hat{f} - f\|_2^2 = \int_0^\infty P(\|\hat{f} - f\|_2^2 > t)\, dt$$
$$\leq \zeta + C_2 \psi_n^2 \int_0^\infty P(\|\hat{f} - f\|_2^2 > C_2 \psi_n^2 t + \zeta)\, dt.$$

With $\tau_n = C_\tau \psi_n^2(1 + t)$, $C_\tau = \xi^{-1}(1 - 2\xi)^2$, this implies that

$$(74) \quad P(\|\hat{f} - f\|_2^2 > C_2 \psi_n^2 t + \zeta) = P(\|\hat{f} - f\|_2^2 > 2\xi \|\hat{f} - f\|_2^2 + \xi \tau_n + \epsilon_n).$$

We have, by Lemma 1, choosing $f^0 = f$, that $\|\hat{f} - f\|_2^2 \leq 2\nu_n[Q(\hat{f} - f)] + \epsilon_n$. This implies that

$$(75) \quad P(\|\hat{f} - f\|_2^2 > C_2 \psi_n^2 t + \zeta) \leq P\left(\sup_{g \in \mathcal{F}} \frac{\nu_n[Q(g - f)]}{w(g)} > \xi\right) \stackrel{\text{def}}{=} P_{\sup},$$

where $w(g) = \|g - f\|_2^2 + \tau_n/2$. We will prove that

$$(76) \quad\quad\quad P_{\sup} \leq \exp(-t \cdot \log_e 2).$$

Together with (73) and (75), this implies the theorem.

*Proof of (76).* We use the peeling device; see, for example, van de Geer (2000), page 69. For $j \geq 0$, let $a_0 = \tau_n/2$, $a_j = 2^{2j} a_0$, $b_j = 2^2 a_j$ and define the following sets of functions: $\mathcal{G}_j = \{g \in \mathcal{F} : a_j \leq w(g) < b_j\}$, $\mathcal{F}_j = \{g \in \mathcal{F} : \|g - f\|_2^2 < b_j\}$. We have that

$$\mathcal{F} = \{g \in \mathcal{F} : w(g) \geq a_0\} = \bigcup_{j=0}^\infty \mathcal{G}_j.$$

Thus,

$$(77) \quad P_{\sup} \leq \sum_{j=0}^\infty P\left(\sup_{g \in \mathcal{G}_j} \frac{\nu_n[Q(g - f)]}{w(g)} > \xi\right)$$
$$\leq \sum_{j=0}^\infty P\left(\sup_{g \in \mathcal{F}_j} \nu_n[Q(g - f)] > \xi a_j\right).$$

By assumption 4 of Theorem 3, $\tilde{G}(\psi_n) = 24\sqrt{2}G(\psi_n)$, where $\tilde{G}$ is defined in (81) for sufficiently large $n$. Thus, by the choice of $C = \xi^{-1}4 \cdot 24\sqrt{2}$ in (24), $\psi_n^2 \geq n^{-1/2}\xi^{-1}4\tilde{G}(\psi_n)$. By the choice of $\xi$, we have that $C_\tau \geq 2$ and thus $a_0 = C_\tau \psi_n^2(1+t)/2 \geq \psi_n^2$. Since $G(\delta)/\delta^2$ is decreasing, by assumption 2 of Theorem 3, $\tilde{G}(\delta)/\delta^2$ is also decreasing and $\xi n^{1/2}/4 \geq \tilde{G}(\psi_n)/\psi_n^2 \geq \tilde{G}(a_0^{1/2})/a_0 \geq \tilde{G}(b_j^{1/2})/b_j$, that is,

$$\xi a_j = \xi b_j/4 \geq n^{-1/2}\tilde{G}(b_j^{1/2}). \tag{78}$$

We now apply Lemma 3 stated below, with (78), to get

$$P\left(\sup_{g \in \mathcal{F}_j} \nu_n[Q(g-f)] > \xi a_j\right) \leq \exp\left\{-\frac{n(\xi a_j)^2 C'}{c^2 b_j^{1-a}}\right\} \tag{79}$$

$$\leq \exp\{-C''(j+1)n\psi_n^{2(1+a)}(1+t)^{1+a}\}, \tag{80}$$

where $C'' = C'c^{-2}\xi^2 2^{2(a-1)}(C_\tau/2)^{1+a}$. Here, we have used the facts that $a_j^2/b_j^{1-a} = 2^{2(a-1)}[2^{2j}C_\tau\psi_n^2(1+t)/2]^{1+a}$ and $2^{2j(a+1)} \geq j+1$. When $0 \leq b \leq 1/2$, we have $\sum_{j=0}^\infty b^{j+1} \leq 2b$. When $n\psi_n^{2(1+a)} \geq (\log_e 2)/C''$, we have $\exp\{-C''n \times \psi_n^{2(1+a)}(1+t)^{1+a}\} \leq 1/2$. Now, we combine (77) and (80) to get the upper bound

$$2\exp\{-C''n\psi_n^{2(1+a)}(1+t)^{1+a}\} \leq \exp\{-t\log_e 2\}.$$

We have proven (76). For the proof of Theorem 3, it remains to prove Lemma 3 below. Lemma 3 gives an exponential tail bound for the Gaussian white noise model.

LEMMA 3. *Let $\nu_n$ be the centered empirical operator of a Gaussian white noise process. Operator $\nu_n$ is defined in (55). Let $\mathcal{G} \subset L_2(\mathbf{R}^d)$ be such that $\sup_{g \in \mathcal{G}} \|g\|_2 \leq R$ and denote by $\mathcal{G}_\delta$ a $\delta$-net of $\mathcal{G}$, $\delta > 0$. Assume that $\delta \mapsto \varrho(Q, \mathcal{G}_\delta)\sqrt{\log_e(\#\mathcal{G}_\delta)}$ is decreasing on $(0, R]$, where $\varrho(Q, \mathcal{G}_\delta)$ is defined in (22) and assume that the entropy integral $G(R)$ defined in (23) is finite. Assume that $\varrho(Q, \mathcal{G}_\delta) = c\delta^{-a}$, where $0 \leq a < 1$ and $c > 0$. Then, for all*

$$\xi \geq n^{-1/2}\tilde{G}(R), \qquad \tilde{G}(R) = \max\{24\sqrt{2}G(R), cR^{1-a}\sqrt{\log_e 2/C'}\} \tag{81}$$

*with $C' = 12^{-2}(C'')^{-2}$ and $C'' = (1-a)^{-3/2}\Gamma(3/2)(\log_e 2)^{-3/2}$, we have*

$$P\left(\sup_{g \in \mathcal{G}} \nu_n(Qg) \geq \xi\right) \leq \exp\left\{-\frac{n\xi^2 C'}{c^2 R^{2-2a}}\right\}.$$

A proof of Lemma 3 is given in the technical report Klemelä and Mammen (2009). The main argument makes use of the chaining technique. An analogous lemma in the direct case is, for example, Lemma 3.2 in van de Geer (2000).

6.5. *Proof of Theorem 4.* The proof is similar to the proof of Theorem 3 up to step (79). At this step, we apply Lemma 4 stated below to get

$$P\Big(\sup_{g\in\mathcal{F}_j}\nu_n[Q(g-f)] > \xi a_j\Big)$$
$$\leq \exp\Big\{-\frac{n(\xi a_j)^2 C'}{c^2 b_j^{1-a}}\Big\} + 2\#\mathcal{G}_{B_2}\exp\Big\{-\frac{1}{12}\frac{n(\xi a_j)^2}{B_\infty c^2 b_j^{1-a} + 2\xi a_j B'_\infty/9}\Big\}.$$

The first term in the right-hand side is handled similarly as in the proof of Theorem 3. For the second term in the right-hand side, we have, for sufficiently large $n$,

$$\exp\Big\{-\frac{1}{12}\frac{n(\xi a_j)^2}{B_\infty c^2 b_j^{1-a} + 2\xi a_j B'_\infty/9}\Big\} = \exp\{-n\psi_n^{2(1+a)}2^{2j}(1+t)^{1+a}C''\}$$

since $a_j^{-a} = (2^{2j}a_0)^{-a} \leq a_0^{-a}$ and $a_0^{-a} \geq 1$ for sufficiently large $n$, and we let $C' = \xi^2 C_\tau^{1+a}/[2^{1+a}12(B_\infty c^2 + 2\xi B'_\infty/9)]$. The proof is completed similarly to the proof of Theorem 3.

We have used Lemma 4, which gives an exponential bound for the tail probability in the case of density estimation.

LEMMA 4. *Let $Y_1,\ldots,Y_n \in \mathbf{R}^d$ be i.i.d. with density $Af$ and let the centered empirical process $\nu_n$ be defined as in (55). Assume that $\|Af\|_\infty \leq B_\infty$. Let $\mathcal{G} \subset L_2(\mathbf{R}^d)$ be such that $\sup_{g\in\mathcal{G}}\|g\|_2 \leq R$. Denote by $\mathcal{G}_\delta$ a $\delta$-bracketing net of $\mathcal{G}$, $\delta > 0$. Let $\mathcal{G}_\delta^L = \{g^L : (g^L, g^U) \in \mathcal{G}_\delta\}$ and $\mathcal{G}_\delta^U = \{g^U : (g^L, g^U) \in \mathcal{G}_\delta\}$. Assume that $\sup_{g\in\mathcal{G}_R^L \cup \mathcal{G}_R^U}\|Qg\|_\infty \leq B'_\infty$. Assume that $\delta \mapsto \varrho_{\mathrm{den}}(Q, \mathcal{G}_\delta)\sqrt{\log_e(\#\mathcal{G}_\delta)}$ is decreasing on $(0, R]$, where $\varrho_{\mathrm{den}}(Q, \mathcal{G}_\delta)$ is defined in (28) and assume that the entropy integral $G(R)$ defined in (29) is finite. Assume that $\varrho_{\mathrm{den}}(Q, \mathcal{G}_\delta) = c\delta^{-a}$, where $0 \leq a < 1$ and $c > 0$. Put $\tilde{G}(R) = B_\infty^{1/2}(9^2 + 96 \cdot 2^{-2a})^{1/2}\max\{24\sqrt{2}G(R), 4 \times (\log_e(2))^{-1}(1-a)^{-3/2}\Gamma(3/2)cR^{1-a}\}$. Then, for all $\xi \geq n^{-1/2}\tilde{G}(R)$, we have*

$$P\Big(\sup_{g\in\mathcal{G}}\nu_n(Qg) \geq \xi\Big)$$
$$\leq 4\exp\Big\{-\frac{n\xi^2 C'}{B_\infty c^2 R^{2-2a}}\Big\}$$
$$+ 2\#\mathcal{G}_R\exp\Big\{-\frac{1}{12}\frac{n\xi^2}{B_\infty c^2 R^{2(1-a)} + 2\xi B'_\infty/9}\Big\},$$

*where $\nu_n$ is the centered empirical process defined in (55).*

A proof of Lemma 4 is given in the technical report Klemelä and Mammen (2009). The proof uses the chaining technique with truncation. The proof follows the techniques developed in Bass (1985), Ossiander (1987), Birgé and Massart (1993), Proposition 3 and van de Geer (2000), Theorem 8.13.

6.6. *Proof of Theorem 5.* We proceed similarly as in the proof of Theorem 1. Choose $f_\delta \in \mathcal{F}_\delta$ such that $\|f - f_\delta\|_2 \leq \delta$, where $f$ is the underlying function in $\mathcal{F}$. Choose $\xi < 1/2$ and put $\zeta = \zeta_1 + \zeta_2$ with $\zeta_1 = (1 - 2\xi)^{-1}(1 + 2\xi)\|f - f_\delta\|_2^2$, $\zeta_2 = \kappa n^{-1} \sum_{j=1}^p \rho_j^2 \lambda_j$ and $\kappa = 4c^{-1}\xi^{-1}(1 - 2\xi)^{-1}$. We have that

$$(82) \qquad E(\|\hat{f} - f\|_2^2) \leq \zeta + \int_0^\infty P(\|\hat{f} - f\|_2^2 > t + \zeta)\, dt.$$

For the integrand of the second term, we have that

$$P(\|\hat{f} - f\|_2^2 > t + \zeta)$$
$$= P(\|\hat{f} - f\|_2^2 > 2\xi\|\hat{f} - f\|_2^2 + (1 - 2\xi)t + (1 - 2\xi)\zeta).$$

We now use Lemma 1. This gives

$$\|\hat{f} - f\|_2^2 \leq \|f - f_\delta\|_2^2 + 2\nu_n(Q(\hat{f} - f_\delta)).$$

Together with the last equalities this gives

$$P(\|\hat{f} - f\|_2^2 > t + \zeta)$$
$$\leq P(\|f - f_\delta\|_2^2 + 2\nu_n(Q(\hat{f} - f_\delta)) > 2\xi\|\hat{f} - f\|_2^2 + (1 - 2\xi)(t + \zeta))$$
$$\leq P(\nu_n(Q(\hat{f} - f_\delta)) > 2^{-1}\xi\|\hat{f} - f_\delta\|_2^2 + 2^{-1}(1 - 2\xi)(t + \zeta_2)).$$

Now, put $w_j = \rho_j / \sum_{l=1}^p \rho_l$ and decompose $f_\delta = f_{\delta,1} + \cdots + f_{\delta,p}$ and $\hat{f} = \hat{f}_1 + \cdots + \hat{f}_p$ with $f_{\delta,j}, \hat{f}_j \in \mathcal{F}_{j,\delta}$. Using assumption (51), we get, with $\beta_j = 2^{-1}(1 - 2\xi)(w_j t + \kappa n^{-1} \rho_j^2 \lambda_j)$, that

$$P(\|\hat{f} - f\|_2^2 > t + \zeta)$$
$$\leq P\left(\sum_{j=1}^p \nu_n(Q(\hat{f}_j - f_{\delta,j})) > 2^{-1}\xi c \sum_{j=1}^p \|\hat{f}_j - f_{\delta,j}\|_2^2 + \sum_{j=1}^p \beta_j\right)$$
$$\leq \sum_{j=1}^p \sum_{g_j \in \mathcal{F}_{j,\delta}} P(\nu_n(Q(g_j - f_{\delta,j})) > 2^{-1}\xi c\|g_j - f_{\delta,j}\|_2^2 + \beta_j).$$

We now use $P(\nu_n(h) > \xi) \leq 2^{-1}\exp(-n\xi^2/[2\|h\|_2^2])$; compare this to the proof of Theorem 1. This gives

$$P(\|\hat{f} - f\|_2^2 > t + \zeta)$$
$$\leq \sum_{j=1}^p \sum_{g_j \in \mathcal{F}_{j,\delta}} 2^{-1} \exp\left[-\frac{n(2^{-1}\xi c\|g_j - f_{\delta,j}\|_2^2 + \beta_j)^2}{2\|Q(g_j - f_{\delta,j})\|^2}\right]$$

$$= \sum_{j=1}^{p} 2^{-1} \exp[-n\xi c 4^{-1}(1-2\xi)w_j \rho_j^{-2} t].$$

By plugging this into (82), we get

$$E(\|\hat{f} - f\|_2^2) \leq \zeta + \sum_{j=1}^{p} \int_0^{\infty} \exp[-n\xi c 4^{-1}(1-2\xi)w_j \rho_j^{-2} t]\, dt$$

$$\leq \zeta + n^{-1} 4[\xi c(1-2\xi)]^{-1} \left(\sum_{j=1}^{p} \rho_j\right)^2.$$

Choosing $\xi = 4^{-1}$ gives the statement of Theorem 5.

**Acknowledgments.** We would like to thank the referees for suggesting improvements and pointing out errors.

Department of Mathematical Sciences  
University of Oulu  
P.O. Box 3000  
90014 University of Oulu  
Finland  
E-mail: klemela@oulu.fi

Department of Economics  
University of Mannheim  
L7 3-5  
68131 Mannheim  
Germany  
E-mail: emammen@rumms.uni-mannheim.de